\documentclass[12pt]{article}
\usepackage{amsfonts,amsmath,amsthm, amssymb}
\usepackage{latexsym, euscript, epic, eepic}
\usepackage{lineno}

\newtheorem{cro}{Corollary}[section]
\newtheorem{defn}{Definition}[section]

\newtheorem{thm}{Theorem}[section]
\newtheorem{lem}{Lemma}[section]

\begin{document}

\title{Multifractal Analysis of Ergodic Averages in Some
Nonuniformly Hyperbolic Systems
 \footnotetext {* Corresponding author}
  }
\author{Xiaoyao Zhou$^1$,  Ercai Chen$^{1,2  *}$ \\
 \small   1 School of Mathematical Sciences and Institute of Mathematics, Nanjing Normal University,\\
  \small   Nanjing 210023, Jiangsu, P.R.China\\
   \small   2 Center of Nonlinear Science, Nanjing University,\\
    \small   Nanjing 210093, Jiangsu, P.R.China\\
     \small    e-mail: zhouxiaoyaodeyouxian@126.com,\\
      \small   ecchen@njnu.edu.cn\\
}
\date{}
\maketitle

\begin{center}
 \begin{minipage}{120mm}
{\small {\bf Abstract.} This article is devoted to the study of the
multifractal analysis of  ergodic averages in some nonuniformly
hyperbolic systems. In particular, our results hold for the robust
classes of multidimensional nonuniformly expanding local
diffeomorphisms and Viana maps.}
\end{minipage}
 \end{center}

\vskip0.5cm
\section{Introduction and Preliminaries}
The multifractal analysis of dynamical systems is a subfield of the
dimension theory of dynamical systems. Roughly speaking,
 multifractal analysis studies the complexity
of the level sets of invariant local quantities obtained from a
dynamical system.

Barreira, Pesin and Schmeling \cite{BarPesSch} introduced the
general concept of multifractal spectrum as follows:

Fix a metric space $X$ and a set $Y$ and let $\varphi:X\to Y$ be a
map. Recently the following problem, known as multifractal analysis
of the map $\varphi,$ has attracted considerable interest. What is
the Hausdorff dimension or the topological entropy or $\ldots$ of
the level sets of $\varphi,$ i.e., What is the Hausdorff dimension
or the topological entropy or $\ldots$ of the following so-called
multifractal decomposition sets of $\varphi$?

\begin{equation}
E(t)=\left\{x\in X:\varphi(x)=t\right\},t\in Y.
\end{equation}
For a topological dynamical system $(X,d,T)$ (or  $(X,T)$ for short)
consisting of a compact metric space $(X,d)$ and a continuous map
$T:X\to X,$ let $
\varphi(x)=\lim\limits_{n\to\infty}\frac{1}{n}\sum\limits_{i=0}^{n-1}\psi(T^ix)$
for some continuous function $\psi:X\to\mathbb{R}.$
  Then there are
fruitful results about the descriptions of the structure (Hausdorff
dimension or topological entropy  or topological pressure) of the
level sets  of $\varphi$ in topological dynamical systems. Early
studies of the level sets was about their dimensions and topological
entropy. See Barreira \& Saussol \cite{BarSau}, Barreira, Saussol \&
Schmeling \cite{BarSauSch}, Oliver \cite{Oli1,Oli2}, Fan \& Feng
\cite{FanFen}, Olsen \cite{Ols2,Ols3,Ols4}, Olsen \& Winter
\cite{OlsWin1, OlsWin2}, Takens \& Verbitskiy \cite{TakVer}, Zhou,
Chen \& Cheng \cite{ZhoCheChe} and Pfister \& Sullivan
\cite{PfiSul2}. Recently, the topological pressures of the level
sets has also been investigated. See Thompson \cite{Tho2}, Pei \&
Chen \cite{PeiChe1}, Yamamoto \cite{Yam1}, Climenhaga \cite{Cli}
 and Zhou \& Chen \cite{ZhoChe1,ZhoChe2}. The reader is referred to
\cite{Bar, Bar2,Bar3} and references therein for recent developments
in multifractal analysis.

Now, nonuniformly hyperbolic systems attract more and more
attentions. We refer the readers to  Barreira \& Pesin
\cite{BarPes}, Chung \&  Takahasi \cite{ChuTak},   Johansson,
Jordan,  Oberg \& Pollicott \cite{JohJorObePol},   Jordan \&   Rams
\cite{JorRam}, Liang, Liao, Sun \& Tian \cite{LiaLiaSUnTia}, Liang,
Liu \& Sun \cite{LiaLiaSUnTia}, Liang, Sun \& Tian \cite{LiaSunTia},
Oliveira \cite{Oliv1}, Oliveira \& Viana \cite{OliVia}, Wang \& Sun
\cite{WanSun} and references therein for recent results in
nonuniformly hyperbolic systems. It is well known that the
specification property plays an important role in some uniformly
hyperbolic dynamical systems. The notion of specification is
slightly weaker than the one introduced by Bowen that requires any
finite sequence of pieces of orbit is well approximated by periodic
orbits. It implies that the dynamical systems have some mixing
property. One should mention that other mild forms of specification
were introduced by
 Pfister \& Sullivan \cite{PfiSul2} and Thompson \cite{Tho3} to the
study of multifractal formalism for Birkhoff averages associated to
beta-shifts, and by Pfister \& Sullivan \cite{PfiSul1}, Yamamoto
\cite{Yam2} and Varandas \cite{Var} to study large deviations.  This
article will use the weak form of specification introduced by
Varandas \cite{Var} in a nonuniformly hyperbolic context.

 Denote by $M(X)$ and $M(X,T)$ the set of all Borel probability
measures on $X$ and the collection of all $T$-invariant Borel
probability measures, respectively. It is well known that $M(X)$ and
$M(X,T)$  equipped with weak* topology are both convex, compact
spaces.

\begin{defn}{\rm\cite{Var}}
We say that $(T,m)$ satisfy the non-uniform specification property
if there exists $\delta>0$ such that for $m$-almost every $x$ and
every $0<\epsilon<\delta$ there exists an integer
$p(x,n,\epsilon)\geq1$ satisfying
\begin{align*}
 \lim\limits_{\epsilon\to0}\limsup\limits_{n\to\infty}\frac{1}{n}p(x,n,\epsilon)=0
  \end{align*}
and so that the following holds: given points $x_1,\cdots, x_k$ in a
full $m$-measure set and positive integers $n_1,\cdots,n_k,$ if
$p_i\geq p(x_i,n_i,\epsilon)$ then there exists $z$ that
$\epsilon$-shadows the orbits of each $x_i$ during $n_i$ iterates
with a time lag of $p(x_i,n_i,\epsilon)$ in between $T^{n_i}(x_i)$
and $x_{i+1},$ that is
\begin{align*}
z\in B_{n_1}(x_1,\epsilon) {\rm~ and~}
T^{n_1+p_1+\cdots+n_{i-1}+p_{i-1}}(z)\in B_{n_i}(x_i,\epsilon)
\end{align*}
for every $2\leq i\leq k,$ where
$B_n(x,\epsilon)=\{y:\max\{d(T^ix,T^iy)\}<\epsilon\}.$
\end{defn}
We may assume that the shadowing property hold on a set $K$ and
$\overline{K}$ is $T$-invariant.

\begin{itemize}
  \item If $m$ is $T$-invariant in the
definition of non-uniform specification, then we let $\overline{K}=$
supp $m.$
  \item If $m$ is Lebesgue measure, then $\overline{K}=X.$
\end{itemize}
So, it is a mild condition that $\overline{K}$ is $T$-invariant.
Now, we state the main result of this article as follows:

\begin{thm}\label{thm}
Suppose $(X,T)$ is a topological dynamical systems with the
non-uniform specification as above. If $\varphi$ and $\psi$ are two
continuous functions on $\overline{K},$ then

\begin{align*}
P(\overline{K}(\varphi,\alpha),\psi)=\sup\left\{h_\nu+\int_{\overline{K}}\psi
d\nu:\nu\in M(\overline{K},T){\rm~and~}\int\varphi d\nu=\alpha\right\},
\end{align*}
where $\overline{K}(\varphi,\alpha):=\left\{x\in
\overline{K}:\frac{1}{n}\sum\limits_{i=0}^{n-1}\varphi(T^ix)=\alpha\right\},
                                     h_\nu$ is the metric entropy of the measure $\nu$ and $P(Z,\psi)$ denotes the topological pressure of $\psi$ with respect
to the set $Z.$
\end{thm}
\section{Proof of Main Result}
In this section, we will prove our main result. The upper bound on
$P(\overline{K}(\varphi,\alpha),\psi)$ is easy to get. However, in
order to obtain the lower bound on
$P(\overline{K}(\varphi,\alpha),\psi),$ the dynamical system should
be endowed with some mixing property such as specification by Takens
\&   Verbitskiy \cite{TakVer}, Tomphson \cite{Tho2}, $g$ almost
product property by Pfister \& Sullivan \cite{PfiSul1,PfiSul2}, Pei
\& Chen \cite{PeiChe1}, Yamamoto \cite{Yam1}. Here, we will use the
weak specification introduced by
 Varandas \cite{Var}. The proof will be divided into the following
two subsection.

Now, we present the definition of topological pressure. Let
$Z\subset X,$ be given and $\Gamma_{n}(Z,\epsilon)$ be the
collection of all finite or countable covers of $Z$ by sets of the
form $B_{m}(x,\epsilon),$ with $m\geq n$. In the next we will denote
$\sum_{i=0}^{n-1}\varphi(T^{i}x)$ by $S_{n}\varphi(x)$.

Set
$$M(Z,t,\varphi,n,\epsilon):=\inf_{\mathcal{C}\in\Gamma_{n}(Z,\epsilon)}\left\{\sum_{B_{m}(x,\epsilon)\in
\mathcal{C}}\exp (-tm+\sup_{y\in
B_{m}(x,\epsilon)}S_{m}\varphi(y))\right\},
$$

and
$$M(Z,t,\varphi,\epsilon)=\lim_{n\to\infty}M(Z,t,\varphi,n,\epsilon),$$

then there exists a unique number $P(Z,\varphi,\epsilon)$ such that
$$P(Z,\varphi,\epsilon)=\inf\{t:M(Z,t,\varphi,\epsilon)=0\}=\sup\{t:M(Z,t,\varphi,\epsilon)=\infty\}.$$
$P(Z,\varphi)=\lim_{\epsilon\to0}P(Z,\varphi,\epsilon)$ is called
the topological pressure of $Z$ with respect to $\varphi$.
\subsection{Upper Bound on $P(\overline{K}(\varphi,\alpha),\psi)$}
$\overline{K}$ can be viewed as a subsystem. The upper bound of
$P(\overline{K}(\varphi,\alpha),\psi)$ holds without extra
assumption.  By \cite{Tho2}, we have

\begin{align*}
P(\overline{K}(\varphi,\alpha),\psi)\leq\sup\left\{h_\nu+\int_{\overline{K}}\psi
d\nu:\nu\in M(\overline{K},T){\rm~and~}\int_{\overline{K}}\varphi
d\nu=\alpha\right\}.
\end{align*}
Before showing the lower bound, we give an important lemma as
follows.

\begin{lem}\cite{Tho1,Tho2}
(Generalised Pressure Distribution Principle) Let $(X,d,T)$ be a
topological dynamical system. Let $Z\subset X$ be an arbitrary Borel
set. Suppose there exist $\epsilon>0$ and $s\geq0$ such that one can
find a sequence of Borel probability measures $\mu_k,$ a constant
$K>0$ and an integer $N$ satisfying

\begin{align*}
\limsup\limits_{k\to\infty}\mu_k(B_n(x,\epsilon))\leq
K\exp(-ns+\sum\limits_{i=0}^{n-1}\psi(T^ix))
\end{align*}
for every ball $B_n(x,\epsilon)$ such that $B_n(x,\epsilon)\cap
Z\neq\emptyset$ and $n\geq N.$ Furthermore, assume that at least one
limit measure $\nu$ of the sequence $\mu_k$ satisfies $\nu(Z)>0.$
Then $P(Z, \psi,\epsilon)\geq s.$

\end{lem}

\subsection{Lower Bound on $P(\overline{K}(\varphi,\alpha),\psi)$}
 The dynamical system needs some mild assumption (non-uniform specification) to obtain the
lower bound of $P(\overline{K}(\varphi,\alpha),\psi).$ Choose a
strictly decreasing sequence $\delta_k\to0$ and fix an arbitrary
$\gamma>0.$ Let us fix $\mu\in M(\overline{K},T)$ satisfying
$\int_{\overline{K}}\varphi d\mu=\alpha$ and
\begin{align*}
h_\mu+\int_{\overline{K}}\psi
d\mu\geq\left\{h_\nu+\int_{\overline{K}}\psi d\nu:\nu\in
M(\overline{K},T){\rm~and~}\int_{\overline{K}}\varphi
d\nu=\alpha\right\}-\gamma.
 \end{align*}
\begin{lem}{\rm \cite{You}\cite{Var}}
For each $\delta_k>0,$ there exists $\eta_k\in M(\overline{K},T)$
such that $\eta_k=\sum\limits_{i=1}^{j(k)}\lambda_i\eta_i^k,$ where
$\sum\limits_{i=1}^{j(k)}\lambda_i=1$ and $\eta_i^k\in
M^e(\overline{K},T)$ satisfying $\left|\int_{\overline{K}}\varphi
d\mu-\int_{\overline{K}}\varphi d\eta_k\right|<\delta_k$ and
$h_{\eta_k}>h_\mu-\delta_k.$
\end{lem}
By Birkhoff's ergodic theorem, we can choose a strictly increasing
sequence $l_k\to\infty$ so that each of the sets
\begin{align*}
Y_{k,i}:=\left\{x\in
\overline{K}:\left|\frac{1}{n}S_n\varphi(x)-\int_{\overline{K}}\varphi
d\eta_i^k\right|<\delta_k {\rm~for~all~}n\geq l_k\right\}
\end{align*}
satisfies $\eta_i^k(Y_{k,i})>1-\gamma$ for every
$k\in\mathbb{N},i\in\{1,\cdots,j(k)\}.$

\begin{lem}{\rm \cite{Tho2}}
Let $\gamma'>0.$ For any sufficiently small $\epsilon>0,$ we can
find a sequence $\widehat{n}_k\to\infty$ with
$[\lambda_i\widehat{n}_k]\geq l_k$ and finite sets $S_{k,i}$ so that
each $S_{k,i}$ is a $([\lambda_i\widehat{n}_k],5\epsilon)$ separated
set for $Y_{k,i}$ and
\begin{align*}
M_{k,i}:=\sum\limits_{x\in
S_{k,i}}\exp\left\{\sum\limits_{i=0}^{n_k-1}\psi(f^ix)\right\}
\end{align*}
satisfies
\begin{align*}
M_{k,i}\geq\exp\left\{[\lambda_i\widehat{n}_k]\left(h_{\eta_i^k}+\int\psi
d\eta_i^k-\frac{4}{j(k)}\gamma'\right)\right\}.
\end{align*}

\end{lem}
We choose $\epsilon$ sufficiently small so that the lemma applies
and $Var(\psi,2\epsilon)<\gamma.$ We fix all the ingredients
provided by the lemma. We now use the non-uniform specification
property to define the set $S_k$ as follows. Let $y_i\in S_{k,i}$
and define $x=x(y_1,\cdots,y_{j(k)})$ to be a choice of point which
belongs to $m$-full measure set $K$ and satisfies

\begin{align*}
d_{[\lambda_i\widehat{n}_k]}(y_l,f^{a_l}x)<\frac{\epsilon}{2^k}
\end{align*}
for all $l\in \{1,\cdots, j(k)\}$ where $a_1=0$ and
$a_l=\sum\limits_{i=1}^{l-1}[\lambda_i\widehat{n}_k]+\max\limits_{x\in
S_{k,i}}p(x,[\lambda_i \widehat{n}_k],\epsilon/2^{k+1})$ for $l\in
\{2,\cdots,j(k)\}.$ In fact, by non-uniform specification, we can
first choose $x^\star\in \overline{K},$ such that

\begin{align*}
d_{[\lambda_i\widehat{n}_k]}(y_l,f^{a_l}x^\star)<\frac{\epsilon}{2^{k+1}}
\end{align*}
for all $l\in \{1,\cdots, j(k)\}$ where $a_1=0$ and
$a_l=\sum\limits_{i=1}^{l-1}[\lambda_i\widehat{n}_k]+\max\limits_{x\in
S_{k,i}}p(x,[\lambda_i\widehat{n}_k],\epsilon/2^{k+1})$ for $l\in
\{2,\cdots,j(k)\}.$ Next, by uniform continuity of $\varphi$ on
$\overline{K},$ we choose $x\in K$ such that
$d_{\widehat{n}_k}(x,x^\star)<\frac{\epsilon}{2^{k+1}}.$ Such $x$ is
what we want, if $x$ is more than one, then we choose one and fix
it. Furthermore, fix $\delta>0,$ for sufficiently small $\epsilon,$
$\widehat{n}_k$ can been chosen so large that

\begin{align*}
\max\limits_{x\in
S_{k,i}}p(x,[\lambda_i\widehat{n}_k],\epsilon/2^{k+1})<\frac{\delta}{2^kj(k)}[\lambda_i]\widehat{n}_k.
\end{align*}

Let $S_k$ be the set of all points constructed in this way. Let
$n_k=\sum\limits_{i=1}^{j(k)}[\lambda_i\widehat{n}_k]+\max\limits_{x\in
S_{k,i}}p(x,[\lambda_i\widehat{n}_k],\epsilon/2^{k+1}).$ Then $n_k$
is the amount of time for which the orbit of the points in $S_k$ has
been prescribed and we have $n_k/\widehat{n}_k\to1.$ We can verify
that $S_k$ is $(n_k,4\epsilon)$ separated and so
$\#S_k=\#S_{k,1}\cdots\#S_{k,j(k)}.$ Let $M_k:=M_{k,1}\cdots
M_{k,j(k)}.$

We assume that $\gamma'$ was chosen to be sufficiently small so that
the following lemma holds.

\begin{lem}\label{lem2.3}
We have
\begin{description}
  \item[(i)] for sufficiently large $k, M_k\geq\exp(h_\mu+\int_{\overline{K}}\psi d\mu-\gamma);$
  \item[(ii)] if $x\in S_k,\left|\frac{1}{n_k}S_{n_k}\varphi(x)-\alpha\right|<\delta_k+Var(\varphi,\epsilon/2^k)+1/k.$
\end{description}
\end{lem}
\subsubsection{Construction of the Intermediate Sets
$\{C_k\}_{k\in\mathbb{N}}$ and $\{L_k\}_{k\in\mathbb{N}}$} In this
subsection, we will construct the intermediate sets
$\{C_k\}_{k\in\mathbb{N}}$ and $\{L_k\}_{k\in\mathbb{N}}.$ First, we
choose a sequence $N_k$ which increases to $\infty$ sufficiently
quickly as follows:

Let $N_1=1, L_1=C_1=S_1. $

Let $N_2\geq 2^{n_1+\max\limits_{x\in
S_1}p(x,n_1,\epsilon/4)+n_3+\max\limits_{x\in
S_3}p(x,n_3,\epsilon/2^4)}.$ Enumerate the points in the sets
$S_2=\{x_i^2:i=1,2,\cdots,\#S_2\}$ and consider the set of words of
length $N_2$ with entries in $\{1,2,\cdots,\#S_2\}.$ Each such word
$\underline{i}=(i_1,\cdots,i_{N_2})$  represents a point in
$S_2^{N_2}.$ Using the non-uniform specification property, we can
choose a point $y^\star:=y^\star(i_1,\cdots,i_{N_2})\in\overline{K}$
which satisfies

\begin{align*}
d_{n_2}(x^2_{i_j},T^{a_j}y^\star)<\frac{\epsilon}{2^{3}}
\end{align*}
for all $j\in\{1,\cdots,N_2\}$ where
$a_j=(j-1)(n_2+\max\limits_{x\in S_2} p(x,n_2,\epsilon/2^{3})).$
Moreover, we can choose $y=y(i_1,\cdots,i_{N_2})\in K$ such that
$d_{n_2}(x^2_{i_j},T^{a_j}y)<\frac{\epsilon}{2^{2}}.$   Collect such
$y$ and define

\begin{align*}
C_2=\left\{y(i_1,\cdots,i_{N_2})\in
K:(i_1,\cdots,i_{N_2})\in\{1,\cdots,\#S_2\}^{N_2}\right\}.
\end{align*}
Let $c_2=N_2n_2+(N_2-1)\max\limits_{x\in S_2
}p(x,n_2,\epsilon/2^{3}).$ We construct $L_{2}$ from $L_1$. Let
$x\in L_1$ and $y\in C_{2}.$ Let $t_1=c_1$ and
$t_{2}=t_1+\max\limits_{x\in L_1}p(x,t_1,\epsilon/2^{3})+c_{2}¡£$
Using non-uniform specification property, we can find a point
$z^\star:=z^\star(x,y)\in\overline{K}$ which satisfies

\begin{align*}
d_{t_1}(x,z^\star)<\frac{\epsilon}{2^{3}} {\rm ~and~}
d_{c_{2}}(y,T^{t_1+\max\limits_{x\in
L_k}p(x,t_1,\epsilon/2^{3})}z^\star)<\frac{\epsilon}{2^{3}}.
\end{align*}
Moreover, we can choose $z=z(x,y)\in K$ such that

\begin{align*}
d_{t_1}(x,z)<\frac{\epsilon}{2^{2}} {\rm ~and~}
d_{c_{2}}(y,T^{t_1+\max\limits_{x\in
L_1}p(x,t_1,\epsilon/2^{3})}z)<\frac{\epsilon}{2^{2}}.
\end{align*}
Collect such $z$ and define $L_{2}=\{z(x,y):x\in L_1,y\in C_{2}\}.$
Generally, we let

\begin{align*}
 N_k\geq 2^{\sum\limits_{i=1}^{k-1}N_in_i+(N_i-1)\max\limits_{x\in
S_i}p(x,n_i,\epsilon/2^{i+1})+\max\limits_{x\in
L_i}p(x,t_i,\epsilon/2^{i+2})+n_{k+1}+\max\limits_{x\in
S_{k+1}}p(x,n_{k+1},\epsilon/2^{k+2})}.
\end{align*}
Enumerate the
points in the sets $S_k$ and write them as follows:

\begin{align*}
S_k=\{x_i^k:i=1,2,\cdots,\#S_k\}.
\end{align*}
Consider the set of words of length $N_k$ with entries in
$\{1,2,\cdots,\#S_k\}.$ Each such word
$\underline{i}=(i_1,\cdots,i_{N_k})$  represents a point in
$S_k^{N_k}.$ Using the non-uniform specification property, we can
choose a point $y^\star:=y^\star(i_1,\cdots,i_{N_k})\in\overline{K}$
which satisfies

\begin{align*}
d_{n_k}(x^k_{i_j},T^{a_j}y^\star)<\frac{\epsilon}{2^{k+1}}
\end{align*}
for all $j\in\{1,\cdots,N_k\}$ where
$a_j=(j-1)(n_k+\max\limits_{x\in S_k} p(x,n_k,\epsilon/2^{k+1})).$
Moreover, we can choose $y=y(i_1,\cdots,i_{N_k})\in K$ such that
$d_{n_k}(x^k_{i_j},T^{a_j}y)<\frac{\epsilon}{2^{k}}.$   Collect such
$y$ and define

\begin{align*}
C_k=\left\{y(i_1,\cdots,i_{N_k})\in
K:(i_1,\cdots,i_{N_k})\in\{1,\cdots,\#S_k\}^{N_k}\right\}.
\end{align*}
Let $c_k=N_kn_k+(N_k-1)\max\limits_{x\in S_k
}p(x,n_k,\epsilon/2^{k+1}).$ Then $c_k$ is the amount of time for
which the orbit of points in $C_k$ has been prescribed.

\begin{lem}
Let $\underline{i}$ and $\underline{j}$ be distinct words in
$\{1,2,\cdots,\#S_k\}^{N_k}.$ Then
$d_{c_k}(y(\underline{i}),y(\underline{j}))>3\epsilon.$
\end{lem}
 We construct $L_{k+1}$
from $L_k$ as follows. Let $x\in L_k$ and $y\in C_{k+1}.$ Let
$t_1=c_1$ and $t_{k+1}=t_k+\max\limits_{x\in
L_k}p(x,t_k,\epsilon/2^{k+2})+c_{k+1}¡£$  Using non-uniform
specification property, we can find a point
$z^\star:=z^\star(x,y)\in\overline{K}$ which satisfies

\begin{align*}
d_{t_k}(x,z^\star)<\frac{\epsilon}{2^{k+2}} {\rm ~and~}
d_{c_{k+1}}(y,T^{t_k+\max\limits_{x\in
L_k}p(x,t_k,\epsilon/2^{k+2})}z^\star)<\frac{\epsilon}{2^{k+2}}.
\end{align*}
Moreover, we can choose $z=z(x,y)\in K$ such that

\begin{align*}
d_{t_k}(x,z)<\frac{\epsilon}{2^{k+1}} {\rm ~and~}
d_{c_{k+1}}(y,T^{t_k+\max\limits_{x\in
L_k}p(x,t_k,\epsilon/2^{k+2})}z)<\frac{\epsilon}{2^{k+1}}.
\end{align*}
Collect such $z$ and define $L_{k+1}=\{z(x,y):x\in L_k,y\in
C_{k+1}\}.$ Note that $t_k$ is the amount of time for which the
orbit of points in $L_k$ has been prescribed.

\begin{lem}
For every $x\in L_k$ and distinct $y_1,y_2\in C_{k+1},$ we have

\begin{align*}
d_{t_k}(z(x,y_1),z(x,y_2))<\frac{\epsilon}{2^k}{\rm~and~}d_{t_{k+1}}(z(x,y_1),z(x,y_2))>2\epsilon.
\end{align*}
Thus $L_k$ is a $(t_k,2\epsilon)$ separated set. In particular, if
$z,z'\in L_k,$ then

\begin{align*}
\overline{B}_{t_k}(z,\frac{\epsilon}{2^k})\cap\overline{B}_{t_k}(z',\frac{\epsilon}{2^k})=\emptyset.
\end{align*}
\end{lem}

\begin{lem}
Let $z=z(x,y)\in L_{k+1}.$ We have

\begin{align*}
\overline{B}_{t_{k+1}}(z,\frac{\epsilon}{2^k})\subset\overline{B}_{t_k}(x,\frac{\epsilon}{2^{k-1}}).
\end{align*}
\end{lem}
\subsubsection{Construction of the Fractal $F$ and a Special Sequence
of Measures $\mu_k$} Let $F_k=\bigcup_{x\in
L_k}\overline{B}_{t_k}(x,\frac{\epsilon}{2^{k-1}}).$ It is obvious
that $F_{k+1}\subset F_k.$ Since we have a decreasing sequence of
compact sets, the intersection $F=\bigcap_kF_k$ is non-empty.
Further, every point $p\in F$ can be uniquely represented by a
sequence $\underline{p}=(\underline{p}_1,\underline{p}_2,\cdots)$
where each
$\underline{p}_i=(\underline{p}^i_1,\cdots,\underline{p}^i_{N_i})\in\{1,2,\cdots,\#S_i\}^{N_i}.$
Each point in $L_k$ can be uniquely represented by a finite word
$(\underline{p}_1,\cdots,\underline{p}_k).$ Let
$y(\underline{p}_i)\in C_i.$ Let
$z_1(\underline{p})=y(\underline{p}_1)$ and proceeding inductively,
let
$z_{i+1}(\underline{p})=z(z_i(\underline{p}),y(\underline{p}_{i+1}))\in
L_{i+1}.$ We can also write $z_i(\underline{p})$ as
$z(\underline{p}_1,\cdots,\underline{p}_i).$ Then define
$p:=\pi\underline{p}$ by

\begin{align*}
p=\bigcap\limits_{i\in\mathbb{N}}\overline{B}_{t_i}(z_i(\underline{p}),\frac{\epsilon}{2^{i-1}}).
\end{align*}
It is clear from our construction that we can uniquely represent
every point in $F$ in this way.

\begin{lem}
Given $z=z(\underline{p}_1,\cdots,\underline{p}_k)\in L_k,$ we have
for all $i\in\{1,\cdots,k\}$ and all $l\in\{1,\cdots,N_i\},$

\begin{align*}
d_{n_i}(x^i_{p^i_l},T^{t_{i-1}+\max\limits_{x\in
L_{i-1}}p(x,t_{i-1},\epsilon/2^{i+1})+(l-1)(n_i+\max\limits_{x\in
S_i}p(x, n_i,\epsilon/2^{i+1}))}z)<2\epsilon.
\end{align*}
\end{lem}

We now define the measures on $F$ which yield the required estimates
for the pressure distribution principle. For each $z\in L_k,$ we
associate a number $\EuScript{L}_k(z)\in (0,\infty).$ Using these
numbers as weights, we define, for each $k,$ an atomic measure
centered on $L_k.$ Precisely, if
$z=z(\underline{p}_1,\cdots,\underline{p}_k),$ we define

\begin{align*}
\EuScript{L}_k(z):=\mathcal{L}(\underline{p}_1)\cdots\mathcal{L}(\underline{p}_k),
\end{align*}
where if
$\underline{p}_i=(p_1^i,\cdots,p^i_{N_i})\in\{1,\cdots,\#S_i\}^{N_i},$
then

\begin{align*}
\mathcal{L}(\underline{p}_i):=\prod\limits_{l=1}^{N_i}\exp
S_{n_i}\psi(x^i_{p_l^i}).
\end{align*}
We define $\nu_k:=\sum\limits_{z\in L_k}\delta_z\EuScript{L}_k(z).$
We normalize $\nu_k$ to obtain a sequence of probability measures
$\mu_k.$ More precisely, we let $\mu_k:=\frac{1}{\kappa_k}\nu_k,$
where $\kappa_k$ is the normalizing constant
$\kappa_k:=\sum\limits_{z\in L_k}\EuScript{L}_k(z)=M_1^{N_1}\cdots
M_k^{N_k}.$ In order to prove the main result of this article, we
present some lemmas. The following lemma is similar to Lemma 5.4 of
\cite{TakVer} so we omit the proof.

\begin{lem}
Suppose $\nu$ is a limit measure of the sequence of probability
measures $\mu_k.$ Then $\nu(F)=1.$
\end{lem}

\begin{lem}
For any $p\in F,$ we have
$\lim\limits_{k\to\infty}\frac{1}{t_k}\sum\limits_{i=0}^{t_k-1}\varphi(T^i(p))=\alpha.$
Thus $F\subset \overline{K}(\varphi,\alpha).$
\end{lem}

\begin{proof}
It relies on (ii) of Lemma \ref{lem2.3}. The proof goes like Lemma
5.3 of \cite{TakVer}.
\end{proof}

Let $\EuScript{B}=B_n(q,\epsilon/2)$ be an arbitrary ball which
intersects $F.$ Let $k$ be the unique number which satisfies
$t_k+\max\limits_{x\in L_k}p(x,t_k,\epsilon/2^{k+2})\leq n<t_{k+1}.$
Let $j\in\{0,\cdots,N_{k+1}-1\}$ be the unique number so

\begin{align*}
t_k+(n_{k+1}+\max\limits_{x\in S_k}p(x,n_k,\epsilon/2^{k+1}))j\leq
n<t_k+(n_{k+1}+\max\limits_{x\in
S_k}p(x,n_k,\epsilon/2^{k+1}))(j+1).
\end{align*}
We assume that $j\geq1$ and  the simpler case $j=0$  is similar.

\begin{lem}\label{lem2.11}
For any $p\geq1,$ suppose $\mu_{k+p}(\EuScript{B})>0,$ there exists
unique $x\in L_k$ and $i_1,\cdots,i_j\in\{1,\cdots,\#S_{k+1}\}$
satisfying

\begin{align*}
&\mu_{k+p}(\EuScript{B})\\
\leq&\frac{1}{\kappa_kM^j_{k+1}}\exp(S_n\psi(q)+2nVar(\psi,2\epsilon)+\|\psi\|(\sum\limits_{i=1}^k(N_i-1)\max\limits_{x\in
S_i}p(x,n_i,\epsilon/2^{i+1})\\&+\sum\limits_{i=1}^k\max\limits_{x\in
L_i}p(x,t_i,\epsilon/2^{i+2})+j\max\limits_{x\in
S_{k+1}}p(x,n_{k+1},\epsilon/2^{k+2})) ).
\end{align*}
\end{lem}

\begin{proof}
Case $p=1.$ Suppose $\mu_{k+1}(\EuScript{B})>0,$ there exists unique
$x\in L_k$ and $i_1,\cdots,i_j\in\{1,\cdots,\#S_{k+1}\}$ satisfying

\begin{align*}
\nu_{k+1}(\EuScript{B})\leq\EuScript{L}_k(x)\prod\limits_{l=1}^j\exp
S_{n_{k+1}}\psi(x_{i_l}^{k+1})M_{k+1}^{N_{k+1}-j},
\end{align*}

Case $p>1.$ Similarly, we have

\begin{align*}
\nu_{k+p}(\EuScript{B})\leq\EuScript{L}_k(x)\prod\limits_{l=1}^j\exp
S_{n_{k+1}}\psi(x_{i_l}^{k+1})M_{k+1}^{N_{k+1}-j}M_{k+2}^{N_{k+2}}\cdots
M_{k+p}^{N_{k+p}}.
\end{align*}
Combining with the fact
\begin{align*}
&\EuScript{L}_k(x)\prod\limits_{l=1}^j\exp
S_{n_{k+1}}\psi({x_{i_l}^{k+1}})\\
&\leq\exp
(S_n\psi(q)+2nVar(\psi,2\epsilon)+\|\psi\|(\sum\limits_{i=1}^k(N_i-1)\max\limits_{x\in
S_i}p(x,n_i,\epsilon/2^{i+1})\\&+\sum\limits_{i=1}^k\max\limits_{x\in
L_i}p(x,t_i,\epsilon/2^{i+2})+j\max\limits_{x\in
S_{k+1}}p(x,n_{k+1},\epsilon/2^{k+2}))),
\end{align*}
the desired result follows.
\end{proof}

Let $C:=h_\mu+\int_{\overline{K}}\varphi d\mu.$

\begin{lem}
For sufficiently large $n,$

\begin{align*}
\limsup\limits_{l\to\infty}\mu_l(B_n(q,\epsilon/2))\leq\exp(-n(C-5\gamma)+S_n\psi(q)).
\end{align*}
\end{lem}
\begin{proof}
For sufficiently large $n, $ we have $\kappa_k
M_{k+1}^j\geq\exp((C-2\gamma)n).$ By Lemma \ref{lem2.11}, the
desired result follows.
\end{proof}

Applying the generalized pressure distribution principle and letting
$\epsilon\to0, \gamma\to0,$ we complete the proof.
\section{Some Applications}
In this section, by the work of Paulo Varandas \cite{Var}, Theorem \ref{thm} can be applied (i) to
Maneville-Pomeau map, (ii) to multidimensional local
diffeomorphisms, and (iii) Viana maps.

{\bf Example 1 Maneville-Pomeau map (Intermittency phenomena)} Given
$\alpha\in (0,1),$ let $T:[0,1]\to[0,1]$ be the $C^{1+\alpha}$
transformation of the interval given by

\begin{align*}
T_\alpha=\left\{
            \begin{array}{ll}
            x(1+2^\alpha x^\alpha), &   0\leq x\leq\frac{1}{2};\\
             2x-1, & \frac{1}{2}<x\leq1.
            \end{array}
          \right.
\end{align*}
Paulo Varandas \cite{Var} proved that when $m$ is SRB measure or the
maximal entropy measure, $([0,1],T,m)$ satisfies non-uniform
specification.

\begin{cro}
If $m$ is SRB measure or the maximal entropy measure for
Maneville-Pomeau map, then
\begin{align*}
&P({\rm
supp}~m(\varphi,\alpha),\psi)\\=&\sup\left\{h_\nu+\int_{{\rm
supp~}m}\psi d\nu:\nu\in M({\rm supp~}m,T){\rm~and~}\int_{{\rm
supp~}m}\varphi d\nu=\alpha\right\}.
\end{align*}
If $\psi>0,$ then

\begin{align*}
&BS({\rm
supp}~m(\varphi,\alpha),\psi)\\=&\sup\left\{h_\nu/\int_{{\rm
supp~}m}\psi d\nu:\nu\in M({\rm supp~}m,T){\rm~and~}\int_{{\rm
supp~}m}\varphi d\nu=\alpha\right\},
\end{align*}
where $BS(Z,\psi)$ is the BS dimension of $Z.$
\end{cro}

{\bf Example 2 multidimensional local diffeomorphisms} Let $T_0$ be
an expanding map in $\mathbb{T}^n$ and take a periodic point $p$ for
$T_0.$ Let $T$ be a $C^1$-local diffeomorphism obtained from $T_0$
by a bifurcation in a small neighborhood $U$ of $p$ in such a way
that:

(1) every point $x\in\mathbb{T}^n$ has some preimage outside $U$;

(2) $\|DT(x)^{-1}\|\leq\sigma^{-1}$ for every $x\in
\mathbb{T}^n\setminus U,$ and $\|DT(x)^{-1}\|\leq L$ for every
$x\in\mathbb{T}^n$ where $\sigma>1$ is large enough or $L>0$ is
sufficiently close to1;

(3) $T$ is topologically exact: for every open set $U$ there is
$N\geq1$ for which $T^N(U)=\mathbb{T}^n$

Paulo Varandas \cite{Var} proved that if $m$ be the unique ergodic
equilibrium state for Holder continuous potential $-\log|det DT|,$
then   $(T,m) $ satisfies non-uniform specification.

\begin{cro}
If $m$ is the unique ergodic equilibrium state for Holder continuous
potential $-\log|det DT| $ in multidimensional local
diffeomorphisms, then
\begin{align*}
&P({\rm
supp}~m(\varphi,\alpha),\psi)\\=&\sup\left\{h_\nu+\int_{{\rm
supp~}m}\psi d\nu:\nu\in M({\rm supp~}m,T){\rm~and~}\int_{{\rm
supp~}m}\varphi d\nu=\alpha\right\}.
\end{align*}
If $\psi>0,$ then

\begin{align*}
&BS({\rm
supp}~m(\varphi,\alpha),\psi)\\=&\sup\left\{h_\nu/\int_{{\rm
supp~}m}\psi d\nu:\nu\in M({\rm supp~}m,T){\rm~and~}\int_{{\rm
supp~}m}\varphi d\nu=\alpha\right\}.
\end{align*}
\end{cro}

 {\bf Example 3  Viana maps}
  Viana maps are obtained as $C^3$ small perturbations of
the skew product $\phi_\alpha$ of the cylinder $S^1\times I$ given
by

\begin{align*}
\phi_\alpha(\theta,x)=(d\theta({\rm mod}1),1-ax^2+\alpha
\cos(2\pi\theta))
\end{align*}
where $d\geq16$ is an integer, $a$ is a Misiurewicz parameter for
the quadratic family, and $\alpha$ is small.

Paulo Varandas \cite{Var} proved that when $m$ is SRB measure or
Lebesgue measure, then   Viana maps satisfy non-uniform
specification.

\begin{cro}
If $m$ is  SRB measure or Lebesgue measure for  a Viana map $T$,
then
\begin{align*}
&P ({\rm
supp}~m(\varphi,\alpha),\psi)\\=&\sup\left\{h_\nu+\int_{{\rm
supp~}m}\psi d\nu:\nu\in M({\rm supp~}m,T){\rm~and~}\int_{{\rm
supp~}m}\varphi d\nu=\alpha\right\}.
\end{align*}
If $\psi>0,$ then

\begin{align*}
&BS({\rm
supp}~m(\varphi,\alpha),\psi)\\=&\sup\left\{h_\nu/\int_{{\rm
supp~}m}\psi d\nu:\nu\in M({\rm supp~}m,T){\rm~and~}\int_{{\rm
supp~}m}\varphi d\nu=\alpha\right\}.
\end{align*}
\end{cro}


\noindent {\bf Acknowledgements.}  The authors want to thank Olsen
and Winter for sharing their articles with us. The research was
supported by the National Natural Science Foundation of China (Grant
No. 11271191) and National Basic Research Program of China (Grant
No. 2013CB834100) and the Foundation for Innovative Program of
Jiangsu Province (Grant No. CXZZ12 0380).


\begin{thebibliography}{50}


\bibitem{Bar} L. Barreira, \emph{Thermodynamic Formalism and Applications
to Dimension Theory}, \emph{Progress in Mathematics} {\bf 294},
Birkh\"{a}user, 2011.

\bibitem{Bar2} L. Barreira, \emph{Ergodic Theory, Hyperbolic  Dynamics and Dimension
Theory}, Springer, 2012.

\bibitem{Bar3} L. Barreira, \emph{Dimension and Recurrence in Hyperbolic Dynamics}, Birkhauser Verlag AG, Basel Boston Berlin, 2008.

\bibitem{BarPes} L. Barreira \& Y. Pesin, Nonuniform hyperbolicity,
Cambridge Univ. Press, Cambridge(2007).


\bibitem{BarPesSch} L. Barreira, Y. Pesin \& J. Schmeling, On a general concept of
multifractality: Multifractal spectrum for dimensions, entropies,
and Lyapunov exponents. Multifractal rigidity, \emph{Chaos} {\bf
7(1)} (1997), 27--38.

\bibitem {BarSau} L. Barreira \& B. Saussol, Variational principles
and mixed multifractal spectra,  \emph{Trans. Amer. Math. Soc}. {\bf
353} (2001), 3919--3944.

\bibitem {BarSauSch} L. Barreira, B. Saussol \& J. Schmeling, Higher-dimensional
multifractal analysis,  \emph{J. Math. Pures Appl}. {\bf 81} (2002),
67--91.

\bibitem{ChuTak} Y. Chung \& H. Takahasi, Multifractal formalism for
Benedicks-Carleson quadratic maps,  \emph{Ergodic Theory Dynamical
Systems} DOI:10.1017/etds.2012.188 (2013).



\bibitem{Cli} V. Climenhaga, Topological pressure of simultaneous level
sets, \emph{Nonlinearity} {\bf 26} (2013), 241--268.

\bibitem {FanFen} A. Fan \& D. Feng, On the distribution of long-term time
averages on symbolic space,  \emph{J. Statist. Phys}. {\bf 99}
(2000), 813--856.


\bibitem{FanLiaPey} A. Fan, L. Liao and J. Peyriere. Generic points
in systems of specification and Banach valued Birkhoff ergodic
average, \emph{Discrete Contin. Dyn. Syst.} {\bf 21} (2008),
1103-1128.

\bibitem{JohJorObePol} A. Johansson, T. Jordan, A. Oberg and M.
Pollicott, Multifractal analysis of non-uniformly hyperbolic
systems, \emph{Israel J. Mathematics} {\bf 177(1)} (2010), 125-144.

\bibitem{JorRam} T. Jordan \& M. Rams,
Multifractal analysis of weak Gibbs measures for non-uniformly
expanding $C^1$ maps,  \emph{Ergodic Theory Dynamical Systems} {\bf
31(1)} (2011), 143-164.


\bibitem{LiaLiaSUnTia} C. Liang, G. Liao, W. Sun and X. Tian,
Saturated sets for nonuniformly hyperbolic systems. Preprint.

\bibitem{LiaLiuSun} C. Liang, G. Liu \& W. Sun, Approximation
properties on invariant measure and Oseledec splitting in
non-uniformly hyperbolic systems, \emph{Trans. Amer. Math. Soc. }
{\bf 361} (2009), 1543-1579.

\bibitem{LiaSunTia} C. Liang, W. Sun \& X. Tian, Ergodic properties
of invariant measures for $C^{1+\alpha}$ non-uniformly hyperbolic
systems, \emph{Ergodic Theory and Dynamical Systems} {\bf 33(2)}
(2013), 560-584.


\bibitem {Oli1} E. Olivier, Analyse multifractale de fonctions continues, \emph{C. R. Acad.
Sci. Paris Sr. I Math}. {\bf 326} (1998), 1171--1174.

\bibitem {Oli2} E. Olivier, Structure multifractale d'une dynamique non expansive
d $\acute{e}$ finie sur unensemble de Cantor, \emph{C. R. Acad. Sci.
Paris Sr. I Math}. {\bf 331} (2000), 605--610.

\bibitem{Oliv1} K. Oliveira, Every expanding measure has the
non-uniform specification property, \emph{Proc. Amer. Math. Soc.}
{\bf 140} (2012), 1309-1320.

\bibitem{OliVia} K. Oliveira \& M. Viana, Thermodynamical formalism
for robust classes of potentials and non-uniformly hyperbolic maps,
\emph{Ergodic Theory and Dynamical Systems} {\bf 28} (2008),
501-533.


\bibitem{Ols2} L. Olsen, Multifractal analysis of divergence points of
deformed measure theoretical Birkhoff averages, \emph{J. Math. Pures
Appl}. {\bf 82} (2003), 1591--1649.

\bibitem{Ols3} L. Olsen, Multifractal analysis of divergence points of
deformed measure theoretical Birkhoff averages. III,
\emph{Aequationes Math}. {\bf 71} (2006) 29--58.

\bibitem{Ols4} L. Olsen, Multifractal analysis of divergence points of
deformed measure theoretical Birkhoff averages. IV: Divergence
points and packing dimension,  \emph{Bull. Sci. Math}. {\bf 132}
(2008), 650--678.

\bibitem {OlsWin1} L. Olsen \& S. Winter, Normal and non-normal points
of self-similar sets and divergence points of self-similar measures,
\emph{J. London Math. Soc}. {\bf 67} (2003), 103--122.

\bibitem {OlsWin2} L. Olsen \& S. Winter, Multifractal analysis of divergence points of deformed
measure theoretical Birkhoff averages. II: Non-linearity, divergence
points and Banach space valued spectra, \emph{Bull. Sci. Math}.
 {\bf 131} (2007), 518--558.

\bibitem{PeiChe1} Y. Pei \& E. Chen, On the variational principle for the
topological pressure for certain non-compact sets,  \emph{Sci. China
Ser. A} {\bf 53(4)} (2010), 1117--1128.


\bibitem{Pes} Y. Pesin, \emph{Dimension Theory in Dynamical Systems},   Contemporary Views and Applications.
Univ. of Chicago Press, 1997.

\bibitem{PfiSul1} C. Pfister \& W. Sullivan, Large deviations estimates for dynamical systems without the
specification property. Applications to the $\beta$-shifts,
\emph{Nonlinearity} {\bf 18} (2005) 237--261.

\bibitem{PfiSul2} C. Pfister \& W. Sullivan, On the topological
entropy of saturated sets,  \emph{Ergodic Theory Dynam. Systems}
{\bf 27} (2007) 929--956.


\bibitem{TakVer} F. Takens \& E. Verbitskiy, On the variational principle for the topological entropy of certain non-compact sets,
\emph{Ergodic Theory Dynam. Systems} {\bf 23(1)} (2003) 317--348.

\bibitem{Tho1} D. Thompson, The irregular set for maps with the specification property has full topological
pressure, \emph{Dynamical Systems: An International Journal} {\bf
25(1)} (2010) 25--51.

\bibitem{Tho2} D. Thompson, A variational principle for topological pressure for certain non-compact
sets,  \emph{J. Lond. Math. Soc}. {\bf 80(3)} (2009) 585--602.

\bibitem{Tho3} D. Thompson, Irregular sets, the beta-transformation
and the almost specification property, \emph{Trans. Amer. Math.
Soc.} to appear.


\bibitem{Var} P. Varandas. Non-uniform specification and large
deviations for weak Gibbs measures. \emph{J. Stat. Phys.} {\bf 146}
(2012), 330-358.

\bibitem{WanSun} Z. Wang \& W. Sun, Lyapunov exponents of hyperbolic
mesures and hyperbolic period orbits, \emph{Trans. Amer. Math. Soc.
} {\bf 362} (2010), 4267-4282.


\bibitem{Yam1} K. Yamamoto, Topological pressure of the set of
generic points for $\mathbb{Z}^{d}$-actions, \emph{Kyushu J. Math}.
{\bf 63} (2009), 191--208.

\bibitem{Yam2} K. Yamamoto, On the weaker forms of the specification
property and their applications, \emph{Proc. Amer. Math. Soc.} {\bf
137} (2009), 3807-3814.

\bibitem{You} L. Young. Large deviations in dynamical systems,  \emph{Tran. Amer. Math. Soc.} {\bf
318} (1990), 525-543.

\bibitem{ZhoChe1} X. Zhou \& E. Chen, Topological pressure of historic set for $\mathbb{Z}^{d}$-actions, \emph{J. Math. Anal. Appl}. {\bf 389} (2012), 394--402.

\bibitem{ZhoChe2} X. Zhou \& E. Chen, Multifractal analysis for the
historic set in topological dynamical systems, \emph{Nonlinearity}
{\bf26} (2013), 1975-1997.


\bibitem{ZhoCheChe} X. Zhou, E. Chen \& W. Cheng, Packing entropy and divergence
points,  \emph{Dynamical Systems: An International Journal} {\bf
27(3)} (2012), 387--402.
\end{thebibliography}
\end{document}